\documentclass[12pt]{amsart}
\usepackage{amscd,amssymb}
\usepackage[graph,frame,poly,arc]{xy}  
\usepackage[plainpages,backref,urlcolor=blue]{hyperref}

\theoremstyle{plain}
\newtheorem{thm}[subsection]{Theorem}
\newtheorem{lem}[subsection]{Lemma}

\theoremstyle{definition}
\newtheorem{rk}[subsection]{Remark}

\newtheorem{ex}[subsection]{Example}

\numberwithin{equation}{section}
\newcommand{\C}{\mathbb{C}}
\newcommand{\PP}{\mathbb{P}}

\DeclareMathOperator{\rank}{rank}

\DeclareMathOperator{\indeg}{indeg}

\begin{document}

\title [The 0-th Fitting ideal of  the Jacobian  module of a plane curve]
{The 0-th Fitting ideal of  the Jacobian  module of a plane curve}

\author[Alexandru Dimca]{Alexandru Dimca$^{1}$}
\address{Universit\'e C\^ ote d'Azur, CNRS, LJAD and INRIA, France and Simion Stoilow Institute of Mathematics,
P.O. Box 1-764, RO-014700 Bucharest, Romania}
\email{dimca@unice.fr}

\author[Gabriel Sticlaru]{Gabriel Sticlaru}
\address{Faculty of Mathematics and Informatics,
Ovidius University
Bd. Mamaia 124, 900527 Constanta,
Romania}
\email{gabrielsticlaru@yahoo.com }

\thanks{$^1$ This work has been partially supported by the French government, through the $\rm UCA^{\rm JEDI}$ Investments in the Future project managed by the National Research Agency (ANR) with the reference number ANR-15-IDEX-01 and by the Romanian Ministry of Research and Innovation, CNCS - UEFISCDI, grant PN-III-P4-ID-PCE-2016-0030, within PNCDI III}



\begin{abstract}  
We describe the 0-th Fitting ideal of  the Jacobian  module of a plane curve in terms of determinants involving the Jacobian syzygies of this curve. This leads to  new characterizations of maximal Tjurina curves,
that is of non free plane curves  whose global Tjurina number equals the upper bound given by A. du Plessis and C.T.C. Wall.
\end{abstract}
 
\maketitle


\section{Introduction} 

Let $S=\C[x,y,z]$ be the polynomial ring in three variables $x,y,z$ with complex coefficients, and let $C:f=0$ be a reduced curve of degree $d$ in the complex projective plane $\PP^2$. 
We denote by $J_f$ the Jacobian ideal of $f$, i.e. the homogeneous ideal in $S$ spanned by the partial derivatives $f_x,f_y,f_z$ of $f$, and  by $M(f)=S/J_f$ the corresponding graded quotient ring, called the Jacobian (or Milnor) algebra of $f$. The global Tjurina number $\tau(C)$ of the curve $C$ is the degree of the Jacobian ideal, namely
$\tau(C)=\dim M(f)_k$, for $k$ large enough.
Consider the graded $S$-module of Jacobian syzygies of $f$, namely
$$AR(f)=\{(a,b,c) \in S^3 \ : \ af_x+bf_y+cf_z=0\}.$$
In this note we consider the following skew-symmetric $S$-bilinear form
\begin{equation}
\label{det1}
\phi: AR(f)^2 \to S, \ \ \phi(r_1,r_2)=\frac{\Delta(E,r_1,r_2)}{f},
\end{equation}
where, for $r_j=(a_j,b_j,c_j) \in AR(f)$ with $j=1,2$, we denote by
$\Delta(E,r_1,r_2)$
the determinant of the $3 \times 3$ matrix $M(E,r_1,r_2)$ which has as first row $x,y,z$, as second row $a_1,b_1,c_1$ and as third row $a_2,b_2,c_2$. It turns out that $\Delta(E,r_1,r_2)$ is divisible by $f$, see \cite{Dmax,dPW}. The image of this morphism $\phi$ is an ideal in $S$, which we denote by $I(C)$ and plan to investigate in this note.

We say that $C:f=0$ is an {\it $m$-syzygy curve} if the graded $S$-module $AR(f)$ is generated by $m$ homogeneous syzygies, say $\rho_1,\rho_2,...,\rho_m$, with $m$ minimal possible, of degrees $d_j=\deg \rho_j$ ordered such that $$1 \leq d_1\leq d_2 \leq ...\leq d_m.$$ 
In fact, the case $d_1=0$ occurs if and only if $C$ is a union of lines through a point, a situation which is not considered in the sequel.
 We call these degrees the {\it exponents} of the curve $C$ and the syzygies $\rho_1,...,\rho_m$ a {\it minimal set of generators } for the module  $AR(f)$. Note that $d_1=mdr(f)$ is the minimal degree of a non trivial Jacobian syzygy in $AR(f)$.  The image of the map
$$AR(f) \to S(d_1+1-d), \ \  \rho \mapsto \phi(\rho_1,\rho)$$
is a shift of the Bourbaki ideal $B(C,\rho_1)\subset S$ associated with the curve $C$ and the minimal degree syzygy $\rho_1$, see \cite{DStJump}. This ideal defines a $0$-dimensional subscheme in $\PP^2$, which plays a key role in understanding the jumping lines of the rank 2 vector bundle $E_C$ associated to the graded $S$-module $AR(f)$.

 For any reduced plane curve $C:f=0$, let  $I_f$ denote the saturation of the ideal $J_f$ with respect to the maximal ideal ${\bf m}=(x,y,z)$ in $S$ and consider the following  local cohomology group, usually called the Jacobian module of $f$, see \cite{Se},
 $$N(f)=I_f/J_f=H^0_{\bf m}(M(f)).$$
We set $n(f)_k=\dim N(f)_k$ for any integer $k$ and also $\nu(C)=\max _j \{n(f)_j\}$, as in \cite{Drcc}. Note that $C$ is free if and only if $\nu(C)=0$, and $C$ is nearly free if and only if $\nu(C)=1$.
We also  consider the invariant 
$$\sigma(C)=\min \{j   : n(f)_j \ne 0\}= \indeg (N(f)).$$
 The self duality of the graded $S$-module $N(f)$, see \cite{HS,Se, SW}, implies that 
\begin{equation}
\label{dual} 
 n(f)_{k} = n(f)_{T-k},
\end{equation}
for any integer $k$, hence
$n(f)_s \ne 0$ exactly for $s=\sigma(C),..., T-\sigma(C)$.

 The minimal possible value for $m$ is two, since the $S$-module $AR(f)$ has rank two for any $f$. The curve $C$ is called {\it free} when $m=2$, since then  $AR(f)$ is a free module, see for such curves \cite{B+,Dmax,DStFD,DStRIMS,Sim1,Sim2,ST,To}. 
 Note that $C:f=0$ is free exactly when $I(C)=S$, a result going back to K. Saito, see \cite{KS}. There is a similar, very general result, saying that an $R$- module $M$ is projective if and only if its first nontrivial Fitting ideal is the whole ring $R$, see \cite[Proposition 20.8]{Eis0}. When $R$ is a positively graded ring, as for instance the polynomial ring $S$, and $
M$ is a graded $R$-module that is projective, then $M$ is a graded free $R$-module, see  \cite[Exercise 4.11]{Eis0}.

The study of 3-syzygy curves was initiated in \cite{DSt3syz}, after the special case of plus-one generated line arrangements was considered by Takuro Abe in \cite{Abe18} and the nearly free curves were introduced in \cite{DStRIMS}.
For such a 3-syzygy curve $C:f=0$, we prove in the second section that
the ideal $I(C)$ is the annihilator ideal $Ann(N(f))$ of the Jacobian module $N(f)$, which is a cyclic $S$-module in this case, and hence we get an isomorphism
$N(f) \simeq S/I(C)$ up-to a shift in degrees, see Theorem \ref{thm1} for a precise statement.

The general case is treated in the third section, where we show that the ideal $I(C)$ is the $0$-th Fitting ideal of the Jacobian module $N(f)$,
and hence in particular, the quotient $S/I(C)$ is an Artinian ring, 
see Theorem \ref{thm2} below. As the proof of this result shows, this gives new information on the second order Jacobian syzygies of the curve $C:f=0$.

In the final section we use Theorem \ref{thm2} and a basic result by 
H. Ito, A. Noma and  M. Ohno on the maximal minors of a matrix with linear form entries, to give new characterizations of maximal Tjurina curves, i.e. non free curves  whose global Tjurina number $\tau(C)$ equals the upper bound given by du Plessis-Wall in \cite{dPW}.

\section{The case of 3-syzygy curves} 

Our first result applies to any reduced curve $C:f=0$

\begin{lem}
\label{lem1}
Let $C:f=0$ be a reduced plane curve and let  $\rho_1,...,\rho_m$ be a  minimal set of generators  for the module  $AR(f)$.
For any pair $j<k$, one has 
$$g_{jk}=\phi(\rho_j,\rho_k) \ne 0.$$
\end{lem}

\proof
Assume $g_{jk}= 0$, which implies the existence of a
relation
$$aE+b\rho_j+c \rho_k=0,$$
with $a,b,c \in S$.
Apply this equality to $f$, looking at a triple $(u,v,w) \in S^3$ as being a derivation $u\partial_x+v\partial_y +w\partial_z$, and use $E(f)=d \cdot f$, $\rho_j(f)=\rho_k(f)=0$ to get $a=0$.
Now in the relation $b\rho_j+c \rho_k=0,$ we can assume that $b$ and $c$ are relatively prime homogeneous polynomials, $c \ne 0$. If $\deg c=0$,
this shows that $\rho_k$ is a multiple of $\rho_j$, hence the initial system of generators for $AR(f)$ is not minimal, a contradiction.
If $\deg c>0$, it follows that $c$ divides all the three components of $\rho_j$,
and hence $c^{-1}\rho_j \in AR(f)$ has a degree $<d_j$, another contradiction with the minimality of the system of generators.
\endproof

For a a 3-syzygy curve $C:f=0$, there are 3 generating syzygies $\rho_1,\rho_2,\rho_3$,  of degrees $1 \leq d_1\leq d_2 \leq  d_3$,
and there is a unique generating second order syzygy involving  $\rho_1,\rho_2,\rho_3$ of the form
\begin{equation}
\label{eqR1} 
R= h_1\rho_1+h_2\rho_2+h_3\rho_3=0,
\end{equation}
with $h_i$ a homogeneous polynomial of degree $\deg h_i=d_j+d_k-d+1$, for any permutation $(i,j,k)$ of $(1,2,3)$, see \cite[Corollary 2.2]{DSt3syz}. A 3-syzygy curve is a plus-one generated curve if and only if $d_1+d_2=d$, see \cite[Theorem 2.3]{DSt3syz}, and a plus-one generated curve is a nearly free curve if and only if one has
in addition $d_2=d_3$.
  For a 3-syzygy curve one has
$$\sigma(C)=3(d-1)-(d_1+d_2+d_3),$$
see \cite[Theorem 3.8]{DSt3syz}. Our first result on the ideal $I(C)$
is the following.

\begin{thm}
\label{thm1}
Let $C:f=0$ be a 3-syzygy curve. With the above notation, one has
$$I(C)=(h_1,h_2,h_3).$$ Moreover, the ideal $I(C)$ is a 0-dimensional complete intersection, and the graded Artinian Gorenstein ring $S/I(C)$ is isomorphic to the shifted graded Jacobian module $N(f)(3(d-1)-(d_1+d_2+d_3))$.
In particular, the Hilbert series $H(N(f);t)$ of the graded Jacobian module $N(f)$ is given by
$$H(N(f);t)=t^{3(d-1)-(d_1+d_2+d_3)}\frac{(t^{d_1+d_2-d+1}-1)(t^{d_1+d_3-d+1}-1)(t^{d_2+d_3-d+1}-1)}{(t-1)^3}.$$

\end{thm}
\proof
It is clear that the ideal $I(C)$ is generated by the elements $g_{ij}$,
for $1\leq i<j \leq 3$. Note that $g_{ij}$ is a nonzero homogeneous polynomial of degree $\deg g_{ij} =d_i+d_j+1-d$.
One clearly has
$$h_2g_{12}+h_3g_{13}=\phi(\rho_1,h_2\rho_2+h_3\rho_3)=\phi(\rho_1,R)=0,$$
in view of the relation \eqref{eqR1} above. On the other hand, it is known that 
$g.c.d.(h_2,h_3)=1$, see \cite[Corollary 2.2]{DSt3syz}, which implies that
there exists a nonzero constant $c \in \C^*$ such that
$g_{12}=ch_3$ and $g_{13}=-ch_2$. Moreover, we have the following equalities
$$g_{23}=h_2^{-1}\phi(h_2\rho_2,\rho_3)=h_2^{-1}\phi(h_2\rho_2+ h_3\rho_3 ,\rho_3)=h_2^{-1}\phi(-h_1\rho_1 ,\rho_3)=ch_1.$$
This shows the following equality of ideals $I(C)=(h_1,h_2,h_3)$. It is known that a free presentation of the Jacobian module $N(f)$, up-to a shift in degrees,  can be obtained by dualizing
the first nontrivial morphism in the minimal free
\begin{equation}
\label{res2}
0 \to S(-d_1-d_2-d_3) \to \oplus_{j=1} ^3S(1-d-d_j)\to S^3(1-d)  \to S
\end{equation}
of the Milnor algebra $M(f)$, see \cite[Proposition 2.1]{DSt3syz} for this resolution, and \cite[Proposition 1.3]{HS} for the claim about the presentation of $N(f)$. The first nontrivial morphism corresponds to the second order syzygy $R$ in \eqref{eqR1} above, hence it is given by
$$ S(-d_1-d_2-d_3) \to \oplus_{j=1} ^3S(1-d-d_j), \ \ h \mapsto (hh_1,hh_2,hh_3).$$
When we dualize it, we get the morphism
$$\oplus_{j=1} ^3S(d+d_j-1) \to S(d_1+d_2+d_3), \ \ (a,b,c) \mapsto ah_1+bh_2+ch_3,$$
which proves our final claims.
\endproof

\begin{rk}
\label{rk1}
In the special case of nearly free curves, the idea of the above proof was  already used in the proof of \cite[Theorem 2.6]{DStRIMS}.
\end{rk}

\section{Fitting ideals and Jacobian syzygies} 

Consider the general form of the minimal graded resolution for the graded $S$-module 
  $M(f)$, the Milnor algebra of a curve $C:f=0$, that is assumed to be not free, namely
\begin{equation}
\label{res2A}
0 \to \oplus_{i=1} ^{m-2}S(-e_i) \to \oplus_{j=1} ^mS(1-d-d_j)\to S^3(1-d)  \to S,
\end{equation}
with $m\geq 3$, $e_1\leq ...\leq e_{m-2}$ and $d_1\leq ...\leq d_m$.
It follows from \cite[Lemma 1.1]{HS} that one has
\begin{equation}
\label{res2B}
e_j=d+d_{j+2}-1+\epsilon_j,
\end{equation}
for $j=1,...,m-2$ and some integers $\epsilon_j\geq 1$. Using  \cite[Formula (13)]{HS}, it follows that one has
\begin{equation}
\label{res2C}
d_1+d_2=d-1+\sum_{i=1} ^{m-2}\epsilon_j.
\end{equation}
In other words, the second order syzygies are generated by $(m-2)$ relations
\begin{equation}
\label{SOS1}
R_i =  h_{i1}\rho_1+h_{i2}\rho_2+ \dots +h_{im}\rho_m=0,
\end{equation}
where $h_{ik}$ is a homogeneous polynomial in $S$ of degree
$$\deg h_{ik}=e_i-d_k+1-d,$$
when $h_{ik} \ne 0$.
Consider the matrix $H=(h_{ik})$, where $1 \leq i \leq m-2$ and $1 \leq k \leq m$, having $(m-2)$ rows and $m$ columns. The square matrix obtained from $H$ by forgetting the $a$-th and the $b$-th columns, for $a<b$,   is denoted by $H_{ab}$ and we consider the signed minor
$$m_{ab}=(-1)^{a+b}\det H_{ab}.$$
Note that $m_{ab}$ is a homogeneous polynomial in $S$ of degree
$$\deg m_{ab}=\sum_{j=1,m-2}e_j-\sum_{k=1,m}d_k +(d_a+d_b)+(m-2)(1-d),$$
when $m_{ab} \ne 0$.
Using the formulas \eqref{res2B} and \eqref{res2C}, we get
\begin{equation}
\label{res2D}
\deg m_{ab}=d_a+d_b+1-d=\deg g_{ab},
\end{equation}
where  we assume that
$m_{ab} \ne 0$.
Using \cite[Proposition 1.3]{HS}, it follows that $H$ is a presentation matrix for the $S$-module $N(f)$. Hence, by definition, the ideal in $S$ generated by all the minors $m_{ab}$ is the 0-th Fitting ideal $Fitt_0(N(f))$ of the Jacobian module $N(f)$. One has
\begin{equation}
\label{Fit1}
Ann(N(f))^{m-2} \subset Fitt_0(N(f)) \subset Ann(N(f)),
\end{equation}
since $N(f)$ is generated by $(m-2)$ elements, see \cite[Proposition 20.7]{Eis0} or \cite[Theorem A2.53]{Eis} for these inclusions in general.
Since $N(f)$ is killed by a high power ${\bf m} ^N$ of the maximal ideal
$${\bf m}=(x,y,z),$$ it follows that the ideal $Fitt_0(N(f))$ is also an $\bf m$-primary ideal. Let $K$ denote the fraction field of the polynomial ring $S$, and note that the matrix $H$ has rank $m-2$ over $K$, since
$Fitt_0(N(f)) \ne 0$.

\begin{thm}
\label{thm2}
Let $C:f=0$ be an $m$-syzygy curve. With the above notation, one has
$I(C)=Fitt_0(N(f))$. In particular, the ideal $I(C)$ is an $\bf m$-primary ideal.
\end{thm}
\proof
Fix an index $k \in [1,m]$, that is a column in the matrix $H$, and denote by $H_k$ the matrix obtained by deleting the $k$-th column in $H$.
First we note that the matrix $H_k$ has also rank $m-2$ over $K$.
Indeed, if we assume $\rank H_k <m-2$, then there is a linear combination of the rows $R'_j$ in $H_k$, let's say
$$\sum_{i=1,m-2} a_iR'_i=0,$$
where $a_i \in S$, not all zero. Since the matrix $H$ has rank $m-2$,
the corresponding linear combination of rows in $H$, namely
$$\sum_{i=1,m-2} a_i(h_{i1},h_{i2}, \dots ,h_{im}) \in S^m,$$
is nonzero, more precisely it is a vector in $S^m$, all of whose coordinates are zero, except the $k$-th coordinate, which is a polynomial $b\ne0$. Now we recall the second order syzygies in \eqref{SOS1}, and get the following obvious contradiction
$$0=\sum_{i=1,m-2} a_iR_i=b\rho_k \ne 0.$$
Therefore, we have indeed $\rank H_k=m-2$.
For any $i \in [1,m-2]$, we consider the following obvious equality
$$0=\phi(\rho_k,R_i)=-\sum_{j=1,k-1} h_{ij}g_{jk}+\sum_{j=k+1,m}h_{ij}g_{kj}.$$
These equations form a linear system of $(m-2)$ equations over the field $K$, 
with $(m-1)$ unknowns, namely the polynomials $\pm g_{jk}$ for $j \ne k$ (the sign is $-$ if and only if $j<k$).
The matrix of this system is the matrix $H_k$, which has  rank $m-2$ as we have seen above. It follows that there is a rational fraction $r_k \in K$ such that
$$\pm g_{jk}=r_k m_{jk},$$
for all $j \ne k$. If we fix a pair $(j,k)$ with $j \ne k$ and apply the above
construction starting with $\rho_j$ instead of $\rho_k$, we deduce that
the corresponding rational function $r_j$ should satisfy $r_j= \pm r_k$.

We have seen that 
$\deg g_{jk}=\deg m_{jk}$, and hence the rational function $r_k$ must be a quotient $p/q$ of two homogeneous polynomials of the same degree in $S$, without common factors. Since $r_j=\pm r_k=\pm p/q$, we deduce that the polynomial $q$ divides all the minors $m_{jk}$.
Since the ideal  $Fitt_0(N(f))$ is an $\bf m$-primary ideal, this is possible only if $q$ is a non zero constant. It follows that $p/q$ is a constant.
This constant is non zero, since $g_{jk}\ne 0$ for any $j<k$ as shown above in Lemma \ref{lem1}.

\endproof

\begin{rk}
\label{rk2}
\noindent (i) Theorem \ref{thm2} implies Theorem \ref{thm1}. Indeed, for a 3-syzygy curve the Jacobian module $N(f)$ is cyclic, that is it is generated by one element. In this case one has an isomorphism
$S/Ann(N(f))=N(f)$ up-to a shift in degrees, and also $Ann(N(f))=Fitt_0(N(f))$ by the inclusions in \eqref{Fit1}.

\noindent (ii) The fact that the ideal $I(C)$ is an $\bf m$-primary ideal has the following geometric interpretation: for any point $p$ on the smooth surface $F: f(x,y,z)-1=0$ in $\C^3$, which is called the Milnor fiber of $f$, there are two syzygies
$\rho_j,\rho_k$  in the given minimal system of generators for $ AR(f)$, such that the values $\rho_j(p), \rho_k(p) \in \C^3$ span the tangent space $T_pF$.

\noindent (iii) As a corollary of Lemma \ref{lem1} and of the proof of Theorem \ref{thm2}, we see that all the $(m-2)$-minors $m_{ab}$ in the matrix $H$, describing the generators of the second order syzygies, are nonzero.
\end{rk}

\begin{ex}
\label{ex2}
Consider the curve $C: f=y^7+x^7+z(x^2+yz)^3=0$. Then the exponents are $d_1=d_2=d_3=d_4=d_5=6$, the global Tjurina number is 
$\tau(C)=12$, and the Hilbert series of the Jacobian module $N(f)$ is given by
$$H(N(f);t)=t^4(3+9t+{13}t^2+{15}t^3+{15}^4+{13}t^5+9t^6+3t^7).$$
The
Hilbert series for $S/I(C)$ is
$$H(S/I(C);t)=1+3t+6t^2+10t^3+15t^4+21t^5+18t^6+6t^7.$$
In particular, the Artinian algebra $S/I(C)$ is not Gorenstein.
\end{ex}

\section{New characterizations of maximal Tjurina curves} 

We have shown in \cite[Theorem 3.1]{DStMax} that a plane curve $C:f=0$  of degree $d$ and with exponents $d_1, \dots,d_m$ such that $d_1 \geq d/2$ satisfies
$$\tau(C)=(d-1)(d-d_1-1)+d_1^2-{ 2d_1-d+2 \choose 2},$$
i.e. the global Tjurina number $\tau(C)$ of $C$ equals the upper bound given by du Plessis-Wall in \cite{dPW}, see also \cite{E}, if and only if one has 
$$d_1=d_2= \cdots = d_m$$
and $m=2d_1-d+3$. For simplicity, we call such a curve a maximal Tjurina curve of type $(d,d_1)$. We have the following result.

\begin{thm}
\label{thm3}
Let $C:f=0$ be an $m$-syzygy curve, with $m\geq 3$. Then the following properties are equivalent.
\begin{enumerate}

\item $C$ is a maximal Tjurina curve of type $(d,d_1)$.

\item The $(m-2) \times m$ matrix $H=(h_{ij})$
coming from the second order syzygies of $C$ has as entries
only linear forms $h_{ij} \in S_1$.

\item $I(C)={\bf m}^{2d_1-d+1}$. More precisely,  the maximal minors $m_{ij}$ of the matrix $H$
form a basis of the $\C$-vector space $S_{2d_1-d+1}$ of homogeneous polynomials of degree $2d_1-d+1$ in $x,y,z$.

\item The Hilbert series for $S/I(C)$ is
$$H(S/I(C);t)=\sum_{i=0}^{2d_1-d} {i+2 \choose 2}t^{i}.$$

\end{enumerate}
Moreover, when these properties hold, then $m=2d_1-d+3$.

\end{thm}

\proof
The fact that $(1)$ implies $(2)$ follows from \cite[Theorem 3.1]{DStMax}.  
If we assume $(2)$, then using the formula \eqref{res2D}, it follows that $d_1= \dots=d_m$ and $m=2d_1-d+3$. Next, the ideal $I(C)$ is an $\bf m$-primary ideal by Theorem \ref{thm2}, and when the entries of $H$ are linear forms, we can apply the main result in \cite{INO} and get that $(2)$ implies $(3)$.
Now we show that $(3)$ implies $(1)$.
By Theorem \ref{thm2} we know that 
$$I(C)=Fitt_0(N(f)) \subset {\bf m}^{m-2},$$
since 
$$\deg h_{ij}=d_i+d_j-d+1\geq d_1+d_2-d+1 \geq 1$$
by formula \eqref{res2C}. Hence $(3)$ implies $m-2 \leq 2d_1-d+1$. On the other hand,
the ideal $I(C)$ has ${m \choose 2}$ generators $g_{ij}$, while to generate
${\bf m}^{2d_1-d+1}$, we need at least ${2d_1-d+3 \choose 2}$ elements of degree $2d_1-d+1$. It follows that $m=2d_1-d+3$ and moreover $d_1=\dots=d_m$, otherwise $\deg g_{ij} >2d_1-d+1$
for some $i<j$, which is impossible.
Hence we have an $m$-syzygy curve $C$ such that $d_1=\dots=d_m$
and $m=2d_1-d+3$. Then \cite[Theorem 3.1]{DStMax} shows that $C$ is a maximal Tjurina curve of type $(d,d_1)$.
The claims $(3)$ and $(4)$ are obviously equivalent, so the proof is complete.
\endproof

Here are some examples.
\begin{ex}
\label{ex3}
The minimal possible value for $m$ in Theorem \ref{thm3} is 3, and this happens exactly when $C:f=0$ is a nearly free curve of even degree
$d=2d_1$ and $d_1=d_2=d_3$, see \cite[Section 4.1]{DStMax}.
Then Theorem \ref{thm1} implies that the 1-minors, i.e. the entries of the matrix $H=(h_1,h_2,h_3)$ are three linearly independent linear forms,
so $I(C)={\bf m}$.

When $m=4$, the degree $d=2r-1$ is odd, and we have the following
sequence of examples of such curves 
$$C_d:f_d=(y^3 - x^2z)x^{r-3}y^{r-1}+x^d+y^d=0,$$
for $r \geq 3$, see \cite[Proposition 4.3]{DStMax}.
It follows that $I(C)={\bf m}^2$ in this case.

The maximal possible value for $m$ is $d+1$, and in this case $d_1=d-1$. We have the following examples of maximal Tjurina curves of type $(d,d-1)$. When $d=2p$ is even, one can take
$$C_{2p}: f_{2p}=(x^2-yz)^{p-1}yz+x^{2p}+y^{2p}=0,$$
while for $d=2p+1$ odd, consider the curves
$$C_{2p+1}: f_{2p+1}=(x^2-yz)^{p-1}xyz+x^{2p+1}+y^{2p+1}=0,$$
see \cite[Section 4.7]{DStMax}. There it was checked, using the software SINGULAR, see \cite{Sing}, that these curves are maximal Tjurina for $2 \leq p \leq 20$.
It follows that $I(C)={\bf m}^{d-1}$ in these cases.

\end{ex}

\end{document}